\documentclass[12pt]{article}
\usepackage{amsfonts, amssymb, amsmath, mathrsfs, amsthm, latexsym, graphicx}
\usepackage{fullpage}
\usepackage{indentfirst}
\theoremstyle{plain}
\newtheorem{theorem}{Theorem}[section]
\newtheorem{lemma}[theorem]{Lemma}
\newtheorem{corollary}[theorem]{Corollary}
\newtheorem{proposition}[theorem]{Proposition}

\theoremstyle{remark}
\newtheorem{remark}[theorem]{Remark}
\theoremstyle{definition}
\newtheorem{definition}[theorem]{Definition}
\numberwithin{equation}{section} 
\DeclareMathOperator{\Red}{Red}
\DeclareMathOperator{\Irr}{Irr}
\DeclareMathOperator{\CF}{CF}

\begin{document}
\title{Equivalence of labeled graphs and lattices}
\author{Dr.\ A.\ N.\ Bhavale}
\date{}
\maketitle
\begin{center}
hodmaths@moderncollegepune.edu.in \\ Department of Mathematics, \\ Modern College of Arts, Science and Commerce (Autonomous), \\ 
Shivajinagar, Pune 411005, M.S., India.
\end{center}
\begin{abstract}
In $1973$, Harary and Palmer posed the problem of enumeration of labeled graphs on $n \geq 1$ unisolated vertices and $l \geq 0$ edges. In $1997$, Bender et al.\ obtained a recurrence relation representing the sequence  $A054548$(OEIS) of labeled graphs on $n \geq 0$ unisolated vertices containing $q \geq \frac{n}{2}$ edges. In $2020$, Bhavale and Waphare obtained a recurrence relation representing the sequence of fundamental basic blocks on $n \geq 0$ comparable reducible elements, having nullity 
$l \geq \lfloor \frac{n+1}{2} \rfloor$. In this paper, we prove the equivalence of these two sequences. We also provide an edge labeling for a given vertex labeled finite simple graph.
\end{abstract}
Keywords: {Chain, Lattice, Digraph, Labeled graph.}\\
MSC Classification $2020$: {$06A05, 06A06, 05C20, 05C78$}.

\section{Introduction} \label{sec1}

In $1955$, Harary \cite{bib10} obtained the number of non-isomorphic linear and connected graphs with $p$ points and $k$ lines. In $1956$, Gilbert \cite{bib9} enumerated 
connected labeled graphs. In $1973$, Harary and Palmer \cite{bib11} obtained the number of ways to label a graph, and also posed the problem of enumeration of labeled graphs on $n\geq 1$ unisolated vertices and $l\geq 0$ edges. 
In $1997$, Bender et al.\ \cite{bib2} obtained asymptotic number of labeled graphs on $n$ unisolated vertices with $q$ edges.
Bender et al.\ \cite{bib2} also obtained the recurrence relation representing the number $d(n, q)$ of labeled graphs on $n$ unisolated vertices, containing $q$ edges 
as $q d(n, q) = (N-q+1) d(n, q-1) + n(n-1) d(n-1, q-1) + N d(n-2, q-1)$ with boundary conditions $d(0, 0) = 1$, and for non-zero values of $n$ and $q$, 
$d(0, q) = d(n, 0) = 0$, where $\frac{n}{2} \leq q \leq N = \binom{n}{2}$.
In $2008$, Tauraso \cite{bib13} obtained the (triangular) sequence of edge coverings of the complete graph $K_n$ for $n \geq 1$, which is interestingly the sequence $A054548$ (see Sloane \cite{bib3}). 
In $2020$, Bhavale and Waphare \cite{bib1} introduced the concept of a fundamental basic block. Let $\mathscr{F}_n(l)$ be the set of all non-isomorphic fundamental basic blocks of nullity $l$, containing $n$ comparable reducible elements. Bhavale and Waphare \cite{bib1} obtained the recurrence relation representing the number of all non-isomorphic
fundamental basic blocks on $n \geq 1$ comparable reducible elements, having nullity $l$, as $|\mathscr{F}_{n+1}(l)| = \displaystyle \sum_{k=1}^{n} \sum_{j=0}^{k} \binom{n}{j} \binom{n-j}{k-j}|\mathscr{F}_{n-j}(l-k)|$ with $|\mathscr{F}_{0}(0)|=1, |\mathscr{F}_{1}(l)|=0$, where $\lfloor \frac{n+2}{2} \rfloor \leq l \leq \binom{n+1}{2}$. Bhavale and Waphare \cite{bib1} also observed that for $n \geq 2$, the sequence $(f(n,l))$ where $f(n,l) = |\mathscr{F}_n(l)|$, is equivalent to the sequence $(d(n,l))$. But this equivalence is yet to establish. Moreover, both of these sequences are equivalent to $A054548$. In this paper, we prove the equivalence of the sequence $(f(n,l))$ with the sequence $(d(n,l))$. We also provide an edge labeling which is uniquely determined for a given vertex labeling of a finite simple graph/directed graph.

A {\it{labeled graph}} $G=(V(G),E(G))$ is a finite set of graph vertices $V(G)$ with a set of graph edges $E(G)$ of $2$-subsets of $V(G)$. Two labeled graphs $G$ and $H$ are said to be {\it{isomorphic}} if and only if there is a one-one (or injective) map from $V(G)$ onto (or surjective) $V(H)$ which preserves not only adjacency but also the labeling. The {\it nullity of a graph} $G$ is given by $|E(G)|-|V(G)|+c$, where $c$ is the number of connected components of $G$. 

Let $\leq$ be a partial order relation on a non-empty set $P$, and let $(P,\leq)$ be a poset. A subset $S$ of a poset $(P,\leq)$ is subposet of $P$ if $S$ itself is a poset under the same relation $\leq$ of $P$. Elements $x,y \in P$ are said to be $comparable$, if either $x \leq y$ or $y \leq x$. Elements $x, y \in P$ are said to be {\it incomparable}, denoted by $x||y$, if $x, y$ are not comparable. A {\it chain} (or a {\it linear order}) is a poset in which any two elements are comparable. An element $c \in P$ is a {\it{lower bound}} ({\it{upper bound}} ) of $a, b \in P$ if $c \leq a, c \leq b ~ (a \leq c, b\leq c)$. A {\it{meet}} of $a, b \in P$, denoted by $a \wedge b$, is defined as the greatest lower bound of $a$ and $b$. A {\it{join}} of $a, b \in P$, denoted by $a \vee b$, is defined as the least upper bound of $a$ and $b$. An element $b$ in $P$ {\it{covers}} an element $a$ in $P$ if $a<b$, and there is no element $c$ in $P$ such that $a<c<b$. 
This fact is denoted by $a \prec b$ and called as a {\it covering} or an {\it edge}.
If $a \prec b$ then $a$ is called a {\it{lower cover}} of $b$, and $b$ is called an {\it upper cover} of $a$. The graph on a poset $P$ with edges as coverings is called  {\it{cover graph}} and is denoted by $C(P)$. The number of coverings in a chain is called {\it{length}} of the chain. If $C$ is a chain on $n$ elements, that is, if $|C|=n$ then its length is $n-1$. Bhavale and Waphare \cite{bib1} defined {\it nullity of a poset} $P$, denoted by $\eta(P)$, to be the nullity of its cover graph $C(P)$.

A poset $L$ is a {\it lattice} if both $a \wedge b$ and $a \vee b$, exist in $L$, $\forall ~ a, b \in L$. A sublattice of a lattice $L$ is a subset of $L$ which is a lattice under the same operations $\wedge$ and $\vee$ of $L$.  For $a \leq b$, the interval $[a, b] = \{ x \in L | a \leq x \leq b \}$ is a sublattice of $L$. An element $x$ in a lattice $L$ is {\it join-reducible}({\it meet-reducible}) in $L$, if there exist $y, z \in L$ both distinct from $x$, such that $y \vee z = x (y \wedge z = x)$. An element $x$ in a lattice $L$ is {\it reducible}, if it is either join-reducible or meet-reducible. The set of all reducible elements of $L$ is denoted by $\Red(L)$.
An element $x$ is {\it join-irreducible}({\it meet-irreducible}), if it is not join-reducible(meet-reducible); An element $x$ is {\it doubly irreducible}, if it is both join-irreducible and meet-irreducible. Lattices $L_1$ and $L_2$ are {\it{isomorphic}} (in symbol, $L_1 \cong L_2$), and the map $\phi : L_1 \to L_2$ is an {\it{isomorphism}} if and only if $\phi$ is one-one and onto (or bijective), and $a \leq b$ in $L_1$ if and only if $\phi(a) \leq \phi(b)$ in $L_2$. 

For the other definitions, notation, and terminology see \cite{{bib15},{bib11},{bib14}}. The following result is due to Rival \cite{bib4}.
\begin{lemma}\cite{bib4} \label{rival}
If $L$ is a lattice then $L \setminus A$ is a sublattice of $L$ for every subset $A$ of all doubly irreducible elements of $L$.
\end{lemma}
In $1974$, Rival \cite{bib4} introduced and studied the class of dismantlable lattices.
\begin{definition}\cite{bib4}
\textnormal{A finite lattice of order $n$ is called {\it{dismantlable}} if there exists a chain $L_{1} \subset L_{2} \subset \cdots\subset L_{n}(=L)$ of sublattices of $L$ such that $|L_{i}| = i$, for all $i$}.
\end{definition}
Brucker and G{\'e}ly \cite{bib16} characterized dismantlable lattices as follows.
\begin{theorem}\cite{bib16} \label{BG}
A lattice $L$ is dismantlable lattice if and only if there exists a chain of lattices $L_1 \subset L_2 \subset \cdots \subset L_n = L$ such that $L_1$ is a singleton and $L_{i-1} = L_i \setminus \{x\}$, where $x$ is a doubly irreducible element of $L_i$.
\end{theorem}
Thakare et al.\ \cite{bib5} introduced the concept of an {\it{adjunct operation of lattices}}. Suppose $L_1$ and $L_2$ are two disjoint lattices and $(a, b)$ is a pair of elements in $L_1$ such that $a<b$ and $a \not\prec b$. Define the partial order $\leq$ on $L = L_1 \cup L_2$ with respect to the pair $(a,b)$ as follows: $x \leq y$ in $L$ if $x, y \in L_1$ and $x \leq y$ in $L_1$, or $x, y \in L_2$ and $x \leq y$ in $L_2$, or $x \in L_1$, $y \in L_2$ and $x \leq a$ in $L_1$, or $x \in L_2$, $y \in L_1$ and $b \leq y$ in $L_1$. The pair $(a, b)$ is called an {\it adjunct pair} of $L$. If $L$ is an adjunct of lattices $L_1$ and $L_2$ with an adjunct pair $(a, b)$ then it is denoted by 
$L = L_1 ]_{a}^{b} L_2$.

The following structure theorem is due to Thakare et al.\ \cite{bib5}.
\begin{theorem} \cite{bib5} \label{st}
A finite lattice is dismantlable if and only if it is an adjunct of chains.
\end{theorem}
Suppose a dismantlable lattice $L$ is an adjunct of the chains $C_0, C_1, C_2, \ldots, C_r$, where $C_0$ is a maximal chain of $L$. Then an adjunct representation of $L$ is given by $L = C_0 ]_{a_1}^{b_1} C_1 ]_{a_2}^{b_2} C_2 \cdots $\\$ ]_{a_r}^{b_r} C_r$, where $(a_i, b_i); 1 \leq i \leq r$ are adjunct pairs. Moreover, an adjunct representation of lattice $L$ is unique up to the order in which adjunct pairs occur.

\begin{proposition}\label{adjunctnuliity} 
If $L = L_1 ]_{a}^{b} L_2$ then $\eta(L) = \eta(L_1) + \eta(L_2) + 1$.
\end{proposition}
\begin{proof}
Note that $|E(L)| = |E(L_1)|+|E(L_2)|+2$ and $|V(L)| = |V(L_1)|+|V(L_2)|$. Further, $\eta(L_1) + \eta(L_2) = (|E(L_1)|-|V(L_1)|+1) + (|E(L_2)|-|V(L_2)|+1) = (|E(L_1)|+|E(L_2)|+2) - (|V(L_1)|+|V(L_2)|) = |E(L)| - |V(L)| = \eta(L) - 1$. Thus $\eta(L) = \eta(L_1) + \eta(L_2) + 1$. 
\end{proof}
Using Theorem \ref{st} and Proposition \ref{adjunctnuliity}, we get the following result.
\begin{corollary}\label{r-1}
A dismantlable lattice with $n$ elements has nullity $r$ if and only if it is an adjunct of $r+1$ chains.
\end{corollary}
Now given a vertex labeling of a simple graph on $n$ unisolated vertices, we obtain an edge labeling of that graph in the following section. Using this edge labeling, in Section \ref{sec4}, we obtain an equivalence between the sequences $(f(n,l))$ and $(d(n,l))$.

\section{Edge labeling of graphs} \label{sec2}

A labeling of a graph $G$ is an assignment of labels to either the vertices or the edges or both subject to certain conditions.
The origin of labelings can be attributed to Rosa \cite{bib20}.
In this section, given a vertex labeling of the complete graph $K_n$, we obtain an edge labeling of $K_n$ which is determined uniquely. In this regard, we begin with the fact that dictionary order is a linear order. That is, if $S = \{(i, j) | 1 \leq i < j \leq n\}$, with a dictionary order as follows for 
$(i, j), (i', j') \in S, (i, j) \leq (i', j')$, if $i \leq i'$, or if $i = i'$ and $j \leq j'$, then $(S, \leq)$ is a chain.
For $1 \leq r < n$, define $S_r = \{(r, j) | r < j \leq n\}$.
\begin{proposition}\label{p2}
For $1 \leq r < n$, $S_r $ is a subchain of $S$. Further, $|S_r|=n-r$ and $S_i \cap S_j = \phi, \forall i \neq j$.
\end{proposition}
\begin{proof} 
As $S_r \subseteq S$, $S_r$ is a chain. Now, it follows that $|S_r| = n-r$, since $r < j \leq n$ implies that $r+1 \leq j \leq r + (n-r)$. Also, by the definition of $S_r$,  $S_i \cap S_j = \phi, \forall i \neq j$.
\end{proof}
\begin{proposition}
For $n \geq 2$, $S = S_1 \oplus S_2 \oplus \cdots \oplus S_{n-1}$.
\end{proposition}
\begin{proof}
By Proposition \ref{p2}, suppose $S_r$ is the chain $(r,r+1) \prec (r,r+2) \prec \cdots \prec (r,n)$, where $1 \leq r \leq n-1$. It is sufficient to prove that $(r,n) \prec (r+1, r+2), \forall r$, where $1 \leq r \leq n-2$. Suppose $(r,n) \leq (i,j) \leq (r+1,r+2)$. If $(i,j) = (r,n)$ then we are done. Otherwise, $i \neq r$ or $j \neq n$ which implies that $r<i$, since $(r,n) \leq (i,j)$. Also, $(i,j) \leq (r+1,r+2)$ implies that $i \leq r+1$, or $i=r+1$ and $j \leq r+2$. If $i \leq r+1$ then $i=r+1$ as $r<i$. But then we must have $j \leq r+2$. Now $i = r+1 < j \leq r+2$. Therefore $j=r+2$. Hence $(i,j) = (r+1,r+2)$. Thus, $S = S_1 \oplus S_2 \oplus \cdots \oplus S_{n-1}$.
\end{proof}
\begin{proposition}\label{p4}
For fixed $r$, where $1 \leq r < n$, if the pair $(r,n)$ is $l^{\rm th}$ element of $S$ then $l = rn - \binom{r+1}{2}$.
\end{proposition}
\begin{proof}
Observe that for $1 \leq r < n$, $(r,n)$ is the largest element of the subchain $S_r$ of $S$. So consider a subchain $C^{(r)} = S_1 \oplus S_2 \oplus \cdots \oplus S_{r}$ of $S$. Note that, the pair $(r,n)$ is also the largest element of $C^{(r)}$. Therefore $(r,n)$ is the $l^{\rm th}$ element of $S$, and $l$ is one more than the length of the chain $C^{(r)}$. Therefore $l = (|C^{(r)}|-1) + 1 = |C^{(r)}| = |S_1 \oplus S_2 \oplus \cdots \oplus S_r| = \displaystyle \sum_{i=1}^{r}|S_i| = rn - \binom{r+1}{2}$, since by Proposition \ref{p2}, $|S_i| = n-i$ and $S_i \cap S_j = \phi$, for all $i \neq j$.
\end{proof}
\begin{proposition}\label{p5}
If $(i,j) \in S$ is $k^{\rm th}$ element of $S$ then $1 \leq k = (i-1)n - \binom{i}{2} + j - i \leq N = \binom{n}{2}$.
\end{proposition}
\begin{proof}
Consider the chain $S_i : (i,i+1) \prec (i,i+2) \prec \cdots \prec (i,n)$. Now $j^{\rm th}$ pair in $S_i$ is $(i,i+j)$. This implies that the pair $(i,j)$ in $S_i$ is $(j-i)^{\rm th}$ pair in $S_i$. Now by Proposition \ref{p4}, if a pair $(i-1,n)$ is the $l^{\rm th}$ element of $S$, then $l = (i-1)n - \binom{i}{2}$. Note that $(i-1,n)$ is the largest element of $S_{i-1}$. Therefore if $(i,j)$ is $k^{\rm th}$ element of $S$ then $(i,j) \in S_i$ and $k = l + (j-i) = (i-1)n - \binom{i}{2} + j - i$.
\end{proof}
\begin{proposition}\label{p6}
For $1 \leq k \leq N = \binom{n}{2}$, there exists unique $(i,j) \in S$ such that $k = (i-1)n-\binom{i}{2}+j-i$.
\end{proposition}
\begin{proof}
For if, suppose there exist $(i,j), (i',j') \in S$ such that $(i,j) \neq (i',j')$ but 
\begin{equation} \label{eq1}
k=(i-1)n-\binom{i}{2}+j-i=(i'-1)n-\binom{i'}{2}+j'-i'.
\end{equation} 

If $i=i'$ then by equation \eqref{eq1}, we have $j=j'$, and hence $(i,j) = (i',j')$, which is a contradiction. Therefore $i \neq i'$. Without loss of generality, suppose $i<i'$. Now from equation \eqref{eq1}, $j-j' = (i'-i)(n-1) + (i'-i) (\frac{1-(i'+i)}{2}) = (i'-i)(n-\frac{1}{2}-\frac{i'}{2}-\frac{i}{2}) = (i'-i)(\frac{n-i'}{2} + \frac{n-i-1}{2}) \geq 0$, as $1 \leq i < i' \leq n$. This implies that $j \geq j'$, a contradiction, since there exists $j \geq i+1$ such that $j<j'$, as $i<j$ and $i < i+1 \leq i' < j'$.
Thus $(i,j) \neq (i',j')$ implies that $(i-1)n-\binom{i}{2}+j-i \neq (i'-1)n-\binom{i'}{2}+j'-i'$, which is not possible by equation \eqref{eq1}. Hence $(i,j) = (i',j')$.
\end{proof}
Let $J_N = \{i | 1 \leq i \leq N = \binom{n}{2}\}$, where $n \geq 2$. 
Then the following result follows immediately from Proposition \ref{p6}.
\begin{theorem}\label{c1}
Let $n \geq 2$. Define $f : S \to J_N$ by $f(i,j) = (i-1)n - \binom{i}{2} + j - i, \forall (i,j) \in S$. Then $f$ is bijective, and hence $|S| = N = \binom{n}{2}$.
\end{theorem}
\begin{proof}
Clearly $f$ is a well defined function, since if $(i,j) = (i',j') \in S$ then $i=i'$ and $j=j'$. Therefore $f(i,j) = f(i',j')$. Further by Proposition \ref{p6}, $f$ is clearly surjective. Note that, if $(i,j) \neq (i',j')$ in $S$ then we get $(i-1)n-\binom{i}{2}+j-i \neq (i'-1)n-\binom{i'}{2}+j'-i'$. Hence $f$ is injective.
\end{proof}
Let $K_n$ be a complete graph on $n \geq 2$ vertices, say $\{1, 2, \ldots, n\}$. Let $\overrightarrow{K_n}$ be the digraph associated with $K_n$ such that $1 \leq i < j \leq n$ whenever edge $(i,j) \in E(\overrightarrow{K_n})$.
Then using Theorem \ref{c1}, we get the following result.
\begin{corollary} \label{t1}
Consider a directed graph $\overrightarrow{K_n}$ with vertex labeling $1, 2, \ldots , n$. Then for $i < j$, a directed edge $\overrightarrow{(i,j)} \in E(\overrightarrow{K_n})$ can be labeled as $k = (i-1)n - \binom{i}{2} + j - i \in J_N$. Moreover, for any $k \in J_N$, there exists unique directed edge $\overrightarrow{(i,j)} \in E(\overrightarrow{K_n})$.
\end{corollary}
In Corollary \ref{t1}, we obtained unique edge labeling of $\overrightarrow{K}_n$. Using this kind of edge labeling, we obtain in the following an edge labeling for a directed subgraph (on $n$ unisolated vertices) of $\overrightarrow{K}_n$.
\begin{theorem}\label{c2}
Suppose $\overrightarrow{G}$ is a directed subgraph (on $n$ unisolated vertices) of $\overrightarrow{K}_n$ with vertex labeling $\{v_1, v_2, \ldots, v_n\}$. Then for $i<j$, if $\overrightarrow{(v_i,v_j)} \in E(\overrightarrow{G})$ then $\overrightarrow{(v_i,v_j)}$ can be labeled as $k = (i-1)n-\binom{i}{2}+j-i \in J_N$. Moreover, if $k \in J_N$ is a label of an edge $\overrightarrow{e} \in E(\overrightarrow{G})$ then $\overrightarrow{e} = \overrightarrow{(v_i,v_j)}$.
\end{theorem}
If we remove direction of each edge of a subgraph $\overrightarrow{G}$ of $\overrightarrow{K}_n$, then we get a simple subgraph $G$ (on $n$ unisolated vertices) of $K_n$.
Thus using Theorem \ref{c2}, for a given finite simple graph $G$ on $n$ unisolated vertices, we get an edge labeling for the graph $G$, which is uniquely determined.

In order to prove equivalence of the sequences $(f(n,l))$ and $(d(n,l))$, we need to study the class $\mathscr{F}_n(l)$ firstly. For this sake in Section \ref{sec3}, we introduce the concept of a complete fundamental basic block, and obtain in general an adjunct representation of a fundamental basic block.

\section{Fundamental basic blocks} \label{sec3}

Bhavale and Waphare \cite{bib1} defined the class of \lq\lq RC-lattices\rq\rq as class of all lattices such that each member of the class has all the reducible elements comparable.
\begin{theorem}\cite{bib1}\label{bwcd}
A lattice in which all the reducible elements are comparable is a dismantlable lattice.
\end{theorem}
Thus an RC-lattice is dismantlable. Bhavale and Waphare \cite{bib1} also introduced the concepts of a doubly irreducible element, a basic block, and a fundamental basic block. An element of a poset $P$ is {\it doubly irreducible} in $P$, if it has at most one upper cover and at most one lower cover in $P$. Let $\Irr(P)$ denote the set of all doubly irreducible elements in the poset $P$. 
\begin{definition}\cite{bib1} \label{bb}
A poset $B$ is a {\it basic block}, if it is one element, or $\Irr(B) = \phi$, or removal of a doubly irreducible element from $B$ reduces nullity of $B$ by one. 
\end{definition}
\begin{definition}\cite{bib1} \label{fbb}
An RC-lattice $F$ is said to be a {\it fundamental basic block}, if it is a basic block, and all the adjunct pairs in an adjunct representation of $F$ are distinct.
\end{definition}
\begin{definition} \label{cfbb}
For $n \geq 2$, {\it complete fundamental basic block}, denoted by  $\CF(n)$, is the fundamental basic block on $n$ reducible elements, having nullity $N = \binom{n}{2}$.
\end{definition}
Using Definition \ref{cfbb}, Theorem \ref{st} and Proposition \ref{p6}, we get an adjunct representation of a complete fundamental basic block in the following result.
\begin{theorem}\label{adcfn}
Let $C : u_1\leq u_2\leq \cdots\leq u_n$ be the chain of reducible elements of $\CF(n)$. Then an adjunct representation of $\CF(n)$ is given by $\CF(n) = C_0 ]^{u_2}_{u_1}\{c_1\} ]^{u_3}_{u_1}\{c_2\} \cdots ]^{u_n}_{u_1}\{c_{n-1}\} $\\$ ]^{u_3}_{u_{2}}\{c_{n}\} ]^{u_4}_{u_{2}}\{c_{n+1}\} \cdots ]^{u_n}_{u_{2}}\{c_{2n-3}\} ]^{u_4}_{u_{3}}\{c_{2n-2}\} \cdots ]^{u_j}_{u_{i}}\{c_{k}\} \cdots ]^{u_n}_{u_{n-1}}\{c_{N}\}$, where $N = \binom{n}{2}$, $k = (i-1)n-\binom{i}{2}+j-i$, $1\leq i<j\leq n$, and $C_0$ is the chain $u_1 \prec x_1 \prec u_2 \prec x_2 \prec \cdots \prec u_{n-1} \prec x_{n-1} \prec u_n$ with $x_i || c_k$ for $k = (i-1)n-\binom{i}{2}+1$.
\end{theorem}
\begin{proof}
By Theorem \ref{bwcd}, $\CF(n)$ is a dismantlable lattice. Therefore by Theorem \ref{st} and Corollary \ref{r-1}, $\CF(n) = C_0 ]^{b_1}_{a_1}C_1 ]^{b_2}_{a_2}C_2 \cdots ]^{b_N}_{a_N}C_{N}$, where $C_0$ is a maximal chain containing the chain $C : u_1 \leq u_2 \leq \cdots \leq u_n$, $C_i$ is a chain and $(a_i, b_i)$ is an adjunct pair for each $i, 1 \leq i \leq N = \binom{n}{2}$.  Clearly $|C_k|=1, \forall k, 1 \leq k \leq N$. For if, suppose $|C_k| \geq 2$ for some $k$, where $1 \leq k \leq N$. 
Now removal of a doubly irreducible element from $C_k$ results in removal of an (incident) edge from $C_k$, as $|C_k| \geq 2$.
Therefore there exists $y \in \Irr(C_k)$ such that $\eta(\CF(n) \setminus \{y\}) = \eta(\CF(n))$. 
This is a contradiction to the fact that $\CF(n)$ is a basic block. Therefore suppose $C_k = \{c_k\}$, where $1 \leq k \leq N$. 
Now $a_k, b_k \in U = \{u_1, u_2, \ldots, u_n\}, \forall k, 1 \leq k \leq N$. Also by Proposition \ref{p6}, for each $1 \leq k \leq N$, there exists unique $(i,j) \in S$ such that $k = (i-1)n-\binom{i}{2}+j-i$. Hence for each $c_k$ or an adjunct pair $(a_k, b_k)$, there corresponds a unique pair $(u_i, u_j)$ such that $k = (i-1)n-\binom{i}{2}+j-i$, where $1 \leq i < j \leq n$. Moreover, $a_k = u_i$ and $b_k = u_j$.
Now if $[u_i, u_{i+1}] \cap C_0$ consists of more than one doubly irreducible elements, then again we get a contradiction to the fact that $\CF(n)$ is a basic block. Therefore for each $i$, $1 \leq i \leq n-1$, $[u_i, u_{i+1}] \cap C_0$ consists of exactly one doubly irreducible element, say $x_i$. Thus $C_0$ is the chain $u_1 \prec x_1 \prec u_2 \prec x_2 \prec \cdots \prec u_{n-1} \prec x_{n-1} \prec u_n$.
Moreover $x_i || c_k$ for $k = (i-1)n-\binom{i}{2}+1$, as $u_i \prec c_k \prec u_{i+1}$.
\end{proof}
Let $X = C_0 \setminus C$ and $Y = \{c_k | k \in J_N \}$. 
Then from Theorem \ref{adcfn}, $X = \{x_1, x_2, \ldots , x_{n-1}\}$ and as a set $\CF(n) = C_0 \cup Y$. Moreover, $|\CF(n)|=|C_0|+|Y|=|C|+|X|+|Y|=n+(n-1)+N = 2n-1+\binom{n}{2}$.
Also using Theorem \ref{adcfn}, the total number of coverings or edges in $\CF(n)$ is $(2n-2)+ 2 \binom{n}{2}$. Thus it can be verified that $\eta(\CF(n)) = \binom{n}{2}$.

The following result immediately follows from Theorem \ref{adcfn}.
\begin{proposition} \label{N-1}
If $F$ is a fundamental basic block obtained from $\CF(n)$ by removal of a doubly irreducible element, then $F\in \mathscr{F}_n(N-1)$, where $n>2$ and
$N = \eta(\CF(n)) = \binom{n}{2}$.
\end{proposition}
\begin{proof}
Let $z \in \Irr(\CF(n))$, where $n>2$. Let $L = \CF(n) \setminus \{z\}$. Then it is clear that $\Red(L) = \Red(\CF(n))$. 
By Lemma \ref{rival}, $L$ is a sublattice of $\CF(n)$. By Theorem \ref{adcfn}, either $z=c_k$ where $1\leq k\leq N$ or $z=x_i$ where $1\leq i\leq n-1$. If $z=x_i$ then $L \cong \CF(n) \setminus \{c_{k}\}$ for $k = (i-1)n-\binom{i}{2}+1$.
Therefore without loss of generality, we assume that $z=c_k$, where $1 \leq k \leq N$.
As $k \in J_N$, by Proposition \ref{p6}, there exists unique $(i,j) \in S$ such that $k = (i-1)n-\binom{i}{2}+j-i$, where $1 \leq i < j \leq n$. As $c_k \notin L$, $(u_i, u_j)$ does not remain an adjunct pair in an adjunct representation of $L$. Therefore by Theorem \ref{adcfn}, $L$ is an adjunct of $N$ chains, 
and hence by Corollary \ref{r-1}, $\eta(L) = N-1$.
Note that the elements $u_i$ and $u_j$ remains reducible in $L$, since $\Red(L) = \Red(\CF(n))$.
Now if $j > i+1$ then by Definition \ref{fbb} and using Theorem \ref{adcfn}, $L$ itself is the fundamental basic block on $n$ reducible elements with $\eta(L) = N-1$. 
Also, if $j = i+1$ then by Definition \ref{bb}, $L$ does not remain a basic block, since $\eta(L \setminus \{x_i\}) = \eta(L)$.  But in this case, by Definition \ref{fbb} and using Theorem \ref{adcfn}, $L \setminus \{x_i\}$ becomes the fundamental basic block on $n$ reducible elements, having nullity same as that of $L$. 
Note that by Lemma \ref{rival}, $L \setminus \{x_i\}$ is also a sublattice of $L$, and hence of $\CF(n)$.
Thus removal of a doubly irreducible element from $\CF(n)$ gives rise to the fundamental basic block on $n$ reducible elements, having nullity exactly one less than that of $\CF(n)$.
\end{proof}
Note that, removal of an arbitrary number of doubly irreducible elements from $\CF(n)$ does not guarantee the preservation of all the $n$ reducible elements. However, removal of at most $n-2$ doubly reducible elements from $\CF(n)$ preserves all the $n$ reducible elements, since each reducible element is connected to the remaining all $n-1$ reducible elements via $n-1$ adjunct pairs. 

Now by Lemma \ref{rival} and using the repeated application of the treatment which is used in the proof of Proposition \ref{N-1}, we get the following result.
\begin{corollary}\label{r1}
If $F$ is a fundamental basic block obtained from $\CF(n)$ by removal of $N-l$ doubly irreducible elements such that $\Red(F) = \Red(\CF(n))$ where $n>2$, then 
$F \in \mathscr{F}_n(l)$, where $N = \eta(\CF(n)) = \binom{n}{2}$ and $\lfloor \frac{n+1}{2}\rfloor \leq l \leq N$.
Further, if $l \geq N-n+2$ then $f(n, l) = |\mathscr{F}_n(l)| = \binom{N}{l}$.
\end{corollary}
\begin{proof}
In Proposition \ref{N-1}, we have seen that removal of a doubly irreducible element from $\CF(n)$ gives rise to the fundamental basic block on $n$ reducible elements, having nullity exactly one less than that of $\CF(n)$. Now by Theorem \ref{BG}, removal of more than $n-2$ doubly irreducible elements from $\CF(n)$ gives rise to the fundamental basic block, but may contain $\leq n$ reducible elements. But $\Red(F) = \Red(\CF(n))$. Therefore removal of $N-l$ doubly irreducible elements from $\CF(n)$ (which are precisely the elements of the set $Y$) gives rise to the fundamental basic block, say $F$ (on $n$ reducible elements) having nullity $N - (N-l) = l$. Thus $F\in \mathscr{F}_n(l)$.

Now by Lemma \ref{rival}, $F$ is a sublattice of $\CF(n)$. Hence $\eta(F) = l \leq \eta(\CF(n)) = N = \binom{n}{2}$. 
Also, as far as the lower bound on $l$ is concerned, there are the following two cases. \\
Case i: Suppose $n$ is even, say $n=2k, k \in \mathbb{N}$.
By Definition \ref{fbb}, the multiplicity of each adjunct pair in an adjunct representation of $F$ is one. Moreover, an adjunct pair corresponds to exactly two reducible elements of $F$. Therefore, if we want to cover all the $n$ reducible elements of $F$, then at least $k$ adjunct pairs are required. Therefore $l \geq k = \frac{n}{2} = \lfloor \frac{n+1}{2} \rfloor$. Note that, if we consider $l$ adjunct pairs, where $l<k$ then care of all the $2k$ reducible elements can not be taken.\\
Case ii: Suppose $n$ is odd, say $n=2k+1, k \in \mathbb{N}$.
Again in this case, at least $k$ adjunct pairs are required to cover $n-1$ reducible elements. Therefore one more adjunct pair is required to cover all the $n$ reducible elements of $F$, so that none of the reducible elements is left unassigned. Hence in this case, at least 
$k+1 = \frac{n-1}{2} + 1 = \frac{n+1}{2} = \lfloor \frac{n+1}{2} \rfloor$ reducible elements are required. Thus $l \geq \lfloor \frac{n+1}{2} \rfloor$.

Further, if $l \geq N-n+2$ then $N-l \leq n-2$. But then $\Red(F) = \Red(\CF(n))$, and hence $F \in \mathscr{F}_n(l)$. Therefore in this situation, 
$f(n, l) = \binom{N}{N-l}$, which is the number of ways to choose $N-l$ doubly irreducible elements for the removal, out of $N$ doubly irreducible elements of $\CF(n)$ (or the set $Y$). Thus $f(n, l) = \binom{N}{l}$.
\end{proof}
In Theorem \ref{adcfn}, we have obtained an adjunct representation of $\CF(n)$ which is the fundamental basic block on $n$ reducible elements, having maximum nullity 
$N = \binom{n}{2}$. By Lemma \ref{rival} and by Corollary \ref{r1}, each $F\in \mathscr{F}_n(l)$ where $\lfloor \frac{n+1}{2}\rfloor \leq l \leq N$, is a sublattice of $\CF(n)$, since $\Red(F) = \Red(\CF(n))$. Therefore using Theorem \ref{adcfn} and Corollary \ref{r1}, we get an adjunct representation of $F$ in general, by restricting an adjunct representation of $\CF(n)$ to that of $F$. In fact, we have the following result.
\begin{corollary}\label{adf}
An adjunct representation of $F \in \mathscr{F}_n(l)$ is given by $F = C'_0 ]^{b_1}_{a_1}\{c_{q_1}\} ]^{b_2}_{a_2}\{c_{q_2}\} $\\$ \cdots  ]^{b_l}_{a_l}\{c_{q_l}\}$, where for each $s, 1 \leq s \leq l \leq N = \binom{n}{2}$, $a_s=u_i, b_s=u_j$ with $q_s = (i-1)n-\binom{i}{2}+j-i, 1 \leq i < j \leq n$, ($q_1 < q_2 < \cdots < q_l$), and $C'_0$ is the chain containing chain $C : u_1 \leq u_2 \leq \cdots \leq u_n$ of all the reducible elements of $\CF(n)$ and all those $x_i \in X = C_0 \setminus C$ which satisfy $x_i || c_{q_s}$ in $F$ for $q_s = (i-1)n-\binom{i}{2}+1$.
\end{corollary}
\begin{remark} \label{rem1}
Let $Y_F = Y \cap F$, where $Y = \{c_k | k \in J_N \}$. Then by Corollary \ref{adf}, $Y_F = \{c_{q_s} | 1 \leq s \leq l\}$.
Also by Corollary \ref{r1}, if $|Y \setminus Y_{F_1}| = |Y \setminus Y_{F_2}| = N-l$, then $F_1, F_2 \in \mathscr{F}_n(l)$.
Moreover, $F_1 \cong F_2$ if and only if $Y \setminus Y_{F_1} = Y \setminus Y_{F_2}$.
That is, $F_1 \cong F_2$ if and only if $Y_{F_1} = Y_{F_2}$.
\end{remark}

\section{Equivalence of labeled graphs and lattices} \label{sec4}

Let $\mathscr{D}(n,q)$ be the set of all non-isomorphic labeled graphs on $n$ unisolated vertices, containing $q$ edges. Then $\mathscr{D}(n,q)$ is precisely the set of all non-isomorphic labeled subgraphs (on $n$ unisolated vertices, containing $q$ edges) of $K_{n}$. Clearly $d(n, q) = |\mathscr{D}(n,q)|$.
Let $\overrightarrow{\mathscr{D}}(n,q)$ be the set of all non-isomorphic labeled digraphs on $n$ unisolated vertices $v_1, v_2, \ldots , v_n$ such that for $\overrightarrow{G} \in \overrightarrow{\mathscr{D}}(n,q)$, $\overrightarrow{(v_i, v_j)}$ is a directed edge of $\overrightarrow{G}$ whenever $i<j$. Note that, $\overrightarrow{\mathscr{D}}(n,q)$ is precisely the set of all non-isomorphic labeled directed subgraphs (on $n$ unisolated vertices, containing $q$ directed edges) of $\overrightarrow{K}_n$.
\begin{lemma}\label{l2}
For $n \geq 2$ and for $G \in \mathscr{D}(n,q)$, $\lfloor \frac{n+1}{2}\rfloor \leq q \leq \binom{n}{2}$.
\end{lemma}
\begin{proof} 
Let $G \in \mathscr{D}(n,q)$ and $n \geq 2$. Note that the complete graph $K_n$ is the largest (simple) graph on $n$ unisolated vertices which has $N = \binom{n}{2}$ edges. As $G$ is a subgraph of $K_n$, $q \leq N$. Now as far as lower bound on $q$ is concerned, there are the following two cases. \\
Case i: Suppose $n$ is even, say $n=2k, k \in \mathbb{N}$. Since an edge consists of two vertices, at least $k = \frac{n}{2} = \lfloor \frac{n+1}{2} \rfloor$ edges are required to cover all the vertices of $G$, so that none of the vertex remains isolated.\\
Case ii: Suppose $n$ is odd, say $n=2k+1, k \in \mathbb{N}$. Again, in this case at least $k$ edges of $G$ are required to cover $n-1$ vertices out of $n$. Therefore we require one more edge to cover all the $n$ vertices of $G$, so that none of the vertex remains isolated. Thus, in this case at least 
$k+1 = \frac{n-1}{2} + 1 = \lfloor \frac{n+1}{2} \rfloor$ edges are required.
\end{proof}
The following result follows immediately from Lemma \ref{l2}.
\begin{corollary}
For $n \geq 2$ and for $\overrightarrow{G} \in \overrightarrow{\mathscr{D}}(n,q)$, $\lfloor \frac{n+1}{2}\rfloor \leq q \leq \binom{n}{2}$. 
\end{corollary}

\vspace{0.005 cm}
\begin{center}
\unitlength 1mm 
\linethickness{0.4pt}
\ifx\plotpoint\undefined\newsavebox{\plotpoint}\fi 
\begin{picture}(124.273,56.365)(0,0)
\put(15.511,20.742){\circle{3.041}}
\put(15.707,13.739){\circle{3.041}}
\put(15.509,34.182){\circle{3.041}}
\put(15.443,27.779){\circle{3.041}}
\put(15.51,32.598){\line(0,-1){3.408}}
\put(15.475,26.314){\line(0,-1){4.048}}
\put(15.523,19.277){\line(0,-1){4.048}}
\put(15.649,47.807){\circle{3.041}}
\put(15.845,40.804){\circle{3.041}}
\put(15.581,54.845){\circle{3.041}}
\put(15.613,53.379){\line(0,-1){4.048}}
\put(15.662,46.343){\line(0,-1){4.048}}
\put(15.715,39.414){\line(0,-1){3.68}}
\put(24.704,47.614){\circle{3.041}}
\put(24.704,33.913){\circle{3.041}}
\put(24.527,20.744){\circle{3.041}}
\multiput(17.343,40.802)(.046627329,-.033701863){161}{\line(1,0){.046627329}}
\multiput(24.553,32.477)(-.05222973,-.033648649){148}{\line(-1,0){.05222973}}
\multiput(16.823,27.497)(.050348387,-.03356129){155}{\line(1,0){.050348387}}
\multiput(24.479,19.322)(-.042051429,-.033554286){175}{\line(-1,0){.042051429}}
\multiput(17.303,40.743)(.047408805,.033641509){159}{\line(1,0){.047408805}}
\multiput(16.969,54.796)(.045953488,-.033622093){172}{\line(1,0){.045953488}}
\put(5.573,24.797){\circle{3.041}}
\multiput(5.499,23.47)(.033624031,-.0384806202){258}{\line(0,-1){.0384806202}}
\multiput(5.21,26.169)(.0336654412,.0524448529){272}{\line(0,1){.0524448529}}
\put(4.706,34.918){\circle{3.041}}
\multiput(4.921,33.687)(.0336472727,-.0732545455){275}{\line(0,-1){.0732545455}}
\multiput(4.439,36.483)(.0336851211,.0627024221){289}{\line(0,1){.0627024221}}
\put(4.609,48.702){\circle{3.041}}
\multiput(4.921,50.17)(.070098485,.033590909){132}{\line(1,0){.070098485}}
\multiput(4.439,47.374)(.0337357143,-.070225){280}{\line(0,-1){.070225}}
\put(20.921,27.326){\makebox(0,0)[cc]{$u_2$}}
\put(21.114,11.518){\makebox(0,0)[cc]{$u_1$}}
\put(21.403,40.724){\makebox(0,0)[cc]{$u_3$}}
\put(20.728,56.049){\makebox(0,0)[cc]{$u_4$}}
\put(27.96,20.681){\makebox(0,0)[cc]{$c_1$}}
\put(0,34.766){\makebox(0,0)[cc]{$c_3$}}
\put(28.381,34.241){\makebox(0,0)[cc]{$c_4$}}
\put(0,49.377){\makebox(0,0)[cc]{$c_5$}}
\put(0,23.624){\makebox(0,0)[cc]{$c_2$}}
\put(28.486,47.485){\makebox(0,0)[cc]{$c_6$}}
\put(15.175,5.998){\makebox(0,0)[cc]{$\CF(4)$}}
\put(100.575,44.936){\circle{3.041}}
\put(121.896,44.707){\circle{3.041}}
\put(100.804,28.43){\circle{3.041}}
\put(121.666,27.971){\circle{3.041}}
\put(111.263,44.792){\vector(1,0){.07}}\put(102.15,44.792){\line(1,0){18.226}}
\put(122.037,36.423){\vector(0,-1){.07}}\put(121.98,43.301){\line(0,-1){13.755}}
\put(111.434,36.367){\vector(-4,-3){.07}}\multiput(120.719,43.645)(-.0429837963,-.0336967593){432}{\line(-1,0){.0429837963}}
\put(100.258,36.768){\vector(0,-1){.07}}\put(100.316,43.531){\line(0,-1){13.526}}
\put(99.858,49.147){\makebox(0,0)[cc]{$v_1$}}
\put(121.98,49.147){\makebox(0,0)[cc]{$v_2$}}
\put(121.98,23.815){\makebox(0,0)[cc]{$v_3$}}
\put(100.66,24.159){\makebox(0,0)[cc]{$v_4$}}
\put(113.154,48.282){\makebox(0,0)[cc]{$\overrightarrow{e_1}$}}
\put(125.273,35.507){\makebox(0,0)[cc]{$\overrightarrow{e_4}$}}
\put(97.45,36.653){\makebox(0,0)[cc]{$\overrightarrow{e_3}$}}
\put(109.257,31.036){\makebox(0,0)[cc]{$\overrightarrow{e_5}$}}
\put(77.771,20.659){\circle{3.041}}
\put(77.968,13.657){\circle{3.041}}
\put(77.77,34.1){\circle{3.041}}
\put(77.704,27.697){\circle{3.041}}
\put(77.771,32.516){\line(0,-1){3.408}}
\put(77.735,26.232){\line(0,-1){4.048}}
\put(77.784,19.195){\line(0,-1){4.048}}
\put(78.106,40.722){\circle{3.041}}
\put(77.842,54.763){\circle{3.041}}
\put(77.975,39.332){\line(0,-1){3.68}}
\put(86.964,33.831){\circle{3.041}}
\put(86.787,20.662){\circle{3.041}}
\multiput(79.604,40.72)(.046627329,-.033701863){161}{\line(1,0){.046627329}}
\multiput(86.813,32.395)(-.05222973,-.033648649){148}{\line(-1,0){.05222973}}
\multiput(79.083,27.415)(.050354839,-.03356129){155}{\line(1,0){.050354839}}
\multiput(86.739,19.24)(-.042045714,-.033554286){175}{\line(-1,0){.042045714}}
\put(66.966,34.836){\circle{3.041}}
\multiput(67.181,33.605)(.0336472727,-.0732545455){275}{\line(0,-1){.0732545455}}
\multiput(66.699,36.401)(.0336851211,.0627024221){289}{\line(0,1){.0627024221}}
\put(66.869,48.62){\circle{3.041}}
\multiput(67.181,50.088)(.070098485,.033590909){132}{\line(1,0){.070098485}}
\multiput(66.699,47.292)(.0337357143,-.070225){280}{\line(0,-1){.070225}}
\put(83.181,27.244){\makebox(0,0)[cc]{$u_2$}}
\put(83.374,11.436){\makebox(0,0)[cc]{$u_1$}}
\put(83.663,40.642){\makebox(0,0)[cc]{$u_3$}}
\put(82.989,55.967){\makebox(0,0)[cc]{$u_4$}}
\put(77.879,53.365){\line(0,-1){11.374}}
\put(63.177,34.667){\makebox(0,0)[cc]{$c_3$}}
\put(63.177,49.278){\makebox(0,0)[cc]{$c_5$}}
\put(90.401,19.951){\makebox(0,0)[cc]{$c_1$}}
\put(90.821,33.51){\makebox(0,0)[cc]{$c_4$}}
\put(78.079,6.392){\makebox(0,0)[cc]{$F\in \mathscr{F}_4(4)$}}
\put(34.669,43.741){\circle{3.041}}
\put(55.99,43.512){\circle{3.041}}
\put(34.899,27.235){\circle{3.041}}
\put(55.761,26.776){\circle{3.041}}
\put(45.357,43.597){\vector(1,0){.07}}\put(36.244,43.597){\line(1,0){18.226}}
\put(56.132,35.229){\vector(0,-1){.07}}\put(56.074,42.106){\line(0,-1){13.755}}
\put(45.528,35.171){\vector(-4,-3){.07}}\multiput(54.813,42.45)(-.0429837963,-.0336967593){432}{\line(-1,0){.0429837963}}
\put(34.353,35.573){\vector(0,-1){.07}}\put(34.41,42.336){\line(0,-1){13.526}}
\put(33.952,47.952){\makebox(0,0)[cc]{$v_1$}}
\put(56.074,47.952){\makebox(0,0)[cc]{$v_2$}}
\put(56.074,22.62){\makebox(0,0)[cc]{$v_3$}}
\put(34.754,22.964){\makebox(0,0)[cc]{$v_4$}}
\put(45.08,35.099){\vector(4,-3){.07}}\multiput(35.779,42.81)(.0406157205,-.0336724891){458}{\line(1,0){.0406157205}}
\put(45.273,26.713){\vector(-1,0){.07}}\put(54.285,26.713){\line(-1,0){18.024}}
\put(46.863,5.797){\makebox(0,0)[cc]{$\overrightarrow{K}_4$}}
\put(19.707,19.601){\makebox(0,0)[cc]{$x_1$}}
\put(19.601,33.266){\makebox(0,0)[cc]{$x_2$}}
\put(20.022,46.615){\makebox(0,0)[cc]{$x_3$}}
\put(81.5,19.716){\makebox(0,0)[cc]{$x_1$}}
\put(81.815,33.276){\makebox(0,0)[cc]{$x_2$}}
\put(63.326,.149){\makebox(0,0)[cc]{$Figure\ I$}}
\put(111.21,6.727){\makebox(0,0)[cc]{$\overrightarrow{G_F} \in \overrightarrow{\mathscr{D}}(4,4)$}}
\end{picture}
\end{center}
\vspace{0.005 cm}

Consider the chain $C: u_1 \leq u_2 \leq \cdots \leq u_n$ of all the reducible elements of $\CF(n)$, where $n \geq 2$. Then an interval $[u_i, u_j]$ is a sublattice of $\CF(n)$ for $1 \leq i < j \leq n$. Further, $[u_i, u_j] \cong \CF(m)$, where $m = j-i+1 \leq n$. Let $I = \{I_{ij}=[u_i, u_j] | 1 \leq i < j \leq n\}$. By Theorem \ref{c1}, it clearly follows that $|I| = |S| = \binom{n}{2}$.

By Proposition \ref{p5}, each $I_{ij} \in I$ can be associated with a unique directed edge $\overrightarrow{e}_k \in E(\overrightarrow{K_n})$. Conversely, by Proposition \ref{p6}, each directed edge $\overrightarrow{e}_k \in E(\overrightarrow{K_n})$ can be associated with a unique interval $I_{ij} \in I$. 
 
In fact, using Proposition \ref{p6} and Theorem \ref{adcfn}, we get the following result which provides unique association of $\CF(n)$ with $\overrightarrow{K_n}$.

\begin{theorem}\label{p8}
Let $u_1 \leq u_2 \leq \cdots \leq u_n$ be the chain of all reducible elements of $\CF(n)$, where $n \geq 2$.
Consider the class $\mathscr{A} = \{u_i, I_{ij} = [u_i, u_j] | 1 \leq i \leq n, I_{ij} \in I\}$ associated with $\CF(n)$.
Let $V(\overrightarrow{K_n}) = \{v_1, v_2, \ldots, v_n\}$ and $E(\overrightarrow{K_n}) = \{\overrightarrow{e_k} | k \in J_N \}$. 
Define $\psi : \mathscr{A} \to \overrightarrow{K_n}$ by $\psi(u_i) = v_i$ where $1 \leq i \leq n$, and $\psi(I_{ij}) = \psi([u_i, u_j]) = \overrightarrow{(v_i, v_j)} = \overrightarrow{e_k}$ where $k = (i-1)n-\binom{i}{2}+j-i$, $1 \leq i < j \leq n$. Then $\psi$ is a well-defined bijective map.
\end{theorem}
\begin{proof}
To show that $\psi$ is well-defined, it is sufficient to prove $\psi(I_{ij})=\psi(I_{pq})$ in $E(\overrightarrow{K_n})$ whenever $I_{ij}=I_{pq}$ in $I$. So suppose $I_{ij}=I_{pq}$ in $I$.   Therefore $u_i=u_p$ and $u_j=u_q$, since $[u_i, u_j] = [u_p, u_q]$. Also $1 \leq i < j \leq n$, $1 \leq p < q \leq n$. This implies that $i=p$ and $j=q$. Therefore $v_i=v_p, v_j=v_q$, and $k = (i-1)n-\binom{i}{2}+j-i = (p-1)n-\binom{p}{2}+q-p$. Therefore $\overrightarrow{(v_i, v_j)} = \overrightarrow{(v_p, v_q)} = \overrightarrow{e_k}$.
Thus $\psi(I_{ij}) = \psi(I_{pq})$ in $E(\overrightarrow{K_n})$.

To show that $\psi$ is injective, it is sufficient to prove that $I_{ij} = I_{pq}$ in $I$ whenever $\psi(I_{ij}) = \psi(I_{pq})$ in $E(\overrightarrow{K_n})$. So suppose $\psi(I_{ij}) = \psi(I_{pq})$ in $E(\overrightarrow{K_n})$. But then $\overrightarrow{e_k} = \overrightarrow{(v_i, v_j)} = \overrightarrow{(v_p, v_q)}$, with $k = (i-1)n-\binom{i}{2}+j-i = (p-1)n-\binom{p}{2}+q-p$. Therefore $v_i=v_p$ and $v_j=v_q$. Now by Proposition \ref{p6}, $(i,j) = (p,q)$. This implies that $i=p$ and $j=q$. Hence we have $u_i=u_p$ and $u_j=u_q$. Therefore $I_{ij} = [u_i, u_j] = [u_p, u_q] = I_{pq}$ in $I$.

To show that $\psi$ is surjective, it is sufficient to prove that for each $\overrightarrow{e_k} \in E(\overrightarrow{K_n})$, there exists $I_{ij} \in I$ such that $\psi(I_{ij}) = \overrightarrow{e_k}$. So suppose $\overrightarrow{e_k} \in E(\overrightarrow{K_n})$, where $1 \leq k \leq \binom{n}{2}$. By Proposition \ref{p6}, there exists unique $(i,j)$ such that $k = (i-1)n-\binom{i}{2}+j-i$, $1 \leq i < j \leq n$. Thus there exists a unique pair $(u_i, u_j)$ of reducible elements of $F$, and hence a unique interval $[u_i, u_j] = I_{ij} \in I$ such that $\psi(I_{ij}) = \overrightarrow{e_k}$.
\end{proof}
Consider the chain $u_1 \leq u_2 \leq \cdots \leq u_n$ of all reducible elements of $F \in \mathscr{F}_n(l)$, where $n \geq 2$. Let $U = \{u_i | 1 \leq i \leq n\}$.
Now if $l<N$ then there exists $(i, j) \in S$ such that $(u_i, u_j)$ need not be an adjunct pair in an adjunct representation of $F$.
So let $I_F = \{I_{ij} = [u_i, u_j] | 1 \leq i < j \leq n ~ \text{and} ~ (u_i, u_j) ~ \text{is an adjunct pair in} ~ F\}$, and let $\mathscr{A}_F = U \cup I_F$. Then $\mathscr{A}_F = \{u_i, I_{ij} | 1 \leq i \leq n, I_{ij} \in I_F\}$. By Corollary \ref{r-1} and by Corollary \ref{adf}, $|I_F| = l$.
Note that for every $F \in \mathscr{F}_n(l)$, there is a unique class $\mathscr{A}_F$ associated with $F$.
Clearly $\mathscr{A}_F \subseteq \mathscr{A}$. Therefore, if we restrict the domain of the map $\psi$ (see Theorem \ref{p8}) to $\mathscr{A}_F$, then we get the following result.
\begin{corollary}\label{c5}
For $n \geq 2$ and for $F\in \mathscr{F}_n(l)$, the (restricted) map $\psi_F : {\mathscr{A}_F} \to \overrightarrow{K_n}$ given by $\psi_F(u_i) = v_i$, $1 \leq i \leq n$, and for $I_{ij} \in I_F$, $\psi_F(I_{ij}) = \psi_F([u_i, u_j]) = \overrightarrow{(v_i, v_j)} = \overrightarrow{e_k}$, where $k = (i-1)n-\binom{i}{2}+j-i$, $1 \leq i < j \leq n$. Then $\psi_F$ is a well-defined injective map.
\end{corollary}
Using Corollary \ref{adf} and Corollary \ref{c5}, we get the following result.
\begin{proposition} \label{pr8}
For each $F \in \mathscr{F}_n(l)$ there exists unique $\overrightarrow{G_F} \in \overrightarrow{\mathscr{D}}(n,l)$ such that $\psi_F(\mathscr{A}_F) = \overrightarrow{G_F}$.
\end{proposition}
\begin{proof}
By Corollary \ref{c5}, for each $F \in \mathscr{F}_n(l)$, $\psi_F : {\mathscr{A}_F} \to \overrightarrow{K_n}$ is a well defined injective map, where $\mathscr{A}_F = U \cup I_F$.
Suppose $\psi_F(\mathscr{A}_F) = \overrightarrow{G_F}$. Then $\overrightarrow{G_F}$ is a (directed) subgraph of $\overrightarrow{K_n}$ with $V(\overrightarrow{G_F}) = \psi_F(U)$, $E(\overrightarrow{G_F}) = \psi_F(I_F)$, and $\psi_F : {\mathscr{A}_F} \to \overrightarrow{G_F}$ is a bijection. Clearly $|V(\overrightarrow{G_F})| = n$. Also $|E(\overrightarrow{G_F})| = l$, since by Corollary \ref{adf}, $|I_F| = l$.
Thus $\overrightarrow{G_F} \in \overrightarrow{\mathscr{D}}(n,l)$, since none of the vertex of $\overrightarrow{G_F}$ is isolated.
For if, suppose for some $i$, there is an isolated vertex $v_i$ in $V(\overrightarrow{G_F})$. Therefore $\overrightarrow{(v_i, v_j)} \notin E(\overrightarrow{G_F})$ for any vertex $v_j \in V(\overrightarrow{G_F})$, where $j>i$. 
This implies that $\psi_F^{-1}(\overrightarrow{(v_i, v_j)}) = I_{ij} = [u_i, u_j] \notin I_F$ for any reducible element $u_j$ of $F$, where $j>i$. 
In other words, $(u_i, u_j)$ is not an adjunct pair in an adjunct representation of $F$ for any $u_j$, where $j>i$.
Similarly, we can prove that $(u_j, u_i)$ is not an adjunct pair in an adjunct representation of $F$ for any reducible element $u_j$ of $F$, where $j<i$.
Thus $u_i$ does not remain as a reducible element of $F$.
This is a contradiction to the fact that $F$ is a fundamental basic block on $n$ reducible elements.
Thus for any $F \in \mathscr{F}_n(l)$, there is a unique class $\mathscr{A}_F$, and hence there exists a unique directed graph $\overrightarrow{G_F} \in  \overrightarrow{\mathscr{D}}(n,l)$.
\end{proof}
In particular, if $F \in \mathscr{F}_4(4)$ then $\mathscr{A}_F = U \cup I_F$ with $U = \{u_1, u_2, u_3, u_4\}$ and $I_F = \{I_{12}, I_{14}, I_{23}, I_{24}\}$. Also $\overrightarrow{G_F} \in \overrightarrow{\mathscr{D}}(4, 4)$ with $V(\overrightarrow{G_F}) = \{v_1, v_2, v_3, v_4\}$ and $E(\overrightarrow{G_F}) = \{\overrightarrow{e_1}, \overrightarrow{e_3}, \overrightarrow{e_4}, $\\$ \overrightarrow{e_5}\}$ (see {\it Figure I}).

Now using Corollary \ref{adf} and Proposition \ref{pr8}, we get the following result.
\begin{theorem}\label{c6}
 For $n \geq 2$ and for fixed $l$, $\lfloor \frac{n+1}{2}\rfloor \leq l \leq \binom{n}{2}$, there is a one-to-one correspondence between $\mathscr{F}_n(l)$ and $\overrightarrow{\mathscr{D}}(n,l)$.
\end{theorem}
\begin{proof}
Define $\phi : \mathscr{F}_n(l) \to \overrightarrow{\mathscr{D}}(n,l)$ by $\phi(F) = \overrightarrow{G} ~ ( = \overrightarrow{G_F})$, for all $F \in \mathscr{F}_n(l)$. 
By Proposition \ref{pr8}, $\phi$ is a well defined map.
Now using Theorem \ref{p8}, Corollary \ref{c5} and Proposition \ref{pr8}, for each $F \in \mathscr{F}_n(l)$, $\psi_{F} : \mathscr{A}_F \to \overrightarrow{G}$ is a bijection, where $\psi_F$ is a restriction of the bijective map $\psi : \mathscr{A} \to \overrightarrow{K_n}$. 

Let $\overrightarrow{G} \in \overrightarrow{\mathscr{D}}(n,l)$. Then $\overrightarrow{G}$ is a (directed) subgraph of $\overrightarrow{K_n}$. Suppose $V(\overrightarrow{G}) = V(\overrightarrow{K_n}) = \{v_1, v_2, \ldots , v_n\}$ and 
$E(\overrightarrow{G}) = \{\overrightarrow{e_{q_1}}, \overrightarrow{e_{q_2}}, \ldots , \overrightarrow{e_{q_l}} \} \subseteq E(\overrightarrow{K_n}) = \{\overrightarrow{e_k} | k \in J_N\}$. 
By Theorem \ref{c2}, $\overrightarrow{e_k} = \overrightarrow{(v_i, v_j)}$ where $k = (i-1)n - \binom{i}{2} + j - i$.
Now suppose $U = \psi^{-1} (V(\overrightarrow{G}))$ and $J = \psi^{-1} (E(\overrightarrow{G}))$. Then there exists $F \in \mathscr{F}_n(l)$ with $\mathscr{A}_F = U \cup J$ such that $\psi(\mathscr{A}_F) = \overrightarrow{G}$. Thus there exists $F \in \mathscr{F}_n(l)$ with $\mathscr{A}_F = U \cup J$ such that $\phi(F) = \overrightarrow{G}$.
Therefore $\phi$ is a surjective map.

Let $F_1, F_2 \in \mathscr{F}_n(l)$. Suppose $\phi(F_1) = \overrightarrow{G_1}$, $\phi(F_2) = \overrightarrow{G_2}$, and $\overrightarrow{G_1} \cong \overrightarrow{G_2}$.
By Corollary \ref{adf}, suppose $F_1 = C'_0 ]^{b_1}_{a_1}\{c_{q_1}\} ]^{b_2}_{a_2}\{c_{q_2}\} \cdots ]^{b_l}_{a_l}\{c_{q_l}\}$, where for each $s, 1 \leq s \leq l \leq N = \binom{n}{2}$, $a_s=u_i, b_s=u_j$ with $q_s = (i-1)n-\binom{i}{2}+j-i$, $1 \leq i < j \leq n$, ($q_1 < q_2 < \cdots < q_l$), and $C'_0$ is the chain containing chain $C : u_1 \leq u_2 \leq \cdots \leq u_n$ of all reducible elements of $\CF(n)$, and all those $x_i \in X = C_0 \setminus C$ which satisfy $x_i || c_{q_s}$ in $F_1$ for 
$q_s = (i-1)n-\binom{i}{2}+1$. 
Note that $C_0$ is a maximal chain $u_1 \prec x_1 \prec u_2 \prec x_2 \prec \cdots \prec u_{n-1} \prec x_{n-1} \prec u_n$ of $\CF(n)$.
Similarly by Corollary \ref{adf}, suppose $F_2 = C''_0 ]^{b'_1}_{a'_1}\{c_{r_1}\} ]^{b'_2}_{a'_2}\{c_{r_2}\} \cdots ]^{b'_l}_{a'_l}\{c_{r_l}\}$, where for each $s$, $1 \leq s \leq l \leq N = \binom{n}{2}$, $a'_s=u_{i'}, b'_s=u_{j'}$ with $r_s = (i'-1)n-\binom{i'}{2}+j'-i'$, $1 \leq i' <j' \leq n$, ($r_1 < r_2 < \cdots < r_l$), and $C''_0$ is the chain containing chain $C : u_1 \leq u_2 \leq \cdots \leq u_n$ of all reducible elements of $\CF(n)$, and all those $x_i \in X = C_0 \setminus C$ which satisfy $x_i || c_{r_s}$ in $F_2$ for $r_s = (i-1)n-\binom{i}{2}+1$.

As $\overrightarrow{G_1} \cong \overrightarrow{G_2}$, $\psi_{F_1}(\mathscr{A}_{F_1}) \cong \psi_{F_2}(\mathscr{A}_{F_2})$.
As $\psi$ is bijective, $\psi^{-1}(\psi_{F_1}(\mathscr{A}_{F_1})) = \psi^{-1}(\psi_{F_2}(\mathscr{A}_{F_2}))$.
That is, $(\psi^{-1} \circ \psi_{F_1})(\mathscr{A}_{F_1}) = (\psi^{-1} \circ \psi_{F_2})(\mathscr{A}_{F_2})$. 
Hence $\mathscr{A}_{F_1} = \mathscr{A}_{F_2}$, since $\psi^{-1} \circ \psi_{F_1}$ and $\psi^{-1} \circ \psi_{F_2}$ are identity maps on $\mathscr{A}_{F_1}$ and $\mathscr{A}_{F_2}$ respectively.
As $\mathscr{A}_{F_1} = \mathscr{A}_{F_2}$, we have $U \cup I_{F_1} = U \cup I_{F_2}$. 
This implies that $I_{F_1} = I_{F_2}$, {\it that is}, $\{[a_s, b_s] | 1 \leq s \leq l\} = \{[a'_s, b'_s] | 1 \leq s \leq l\}$.

Now suppose for fixed $s$ where $1 \leq s \leq l$, $[a_s, b_s] \cong [a'_t, b'_t]$ for some $t, 1 \leq t \leq l$. Then $a_s = a'_t$ and $b_s = b'_t$. 
But $a_s \prec c_{q_s} \prec b_s$ and $a'_t\prec c_{r_t} \prec b'_t$. Therefore $c_{q_s} = c_{r_t}$, and hence $q_s = r_t$. 
Thus for all $s, 1 \leq s \leq l$, there exists $t$, where $1 \leq t \leq l$ such that $c_{q_s} = c_{r_t}$. Therefore $Y_{F_1} = Y \cap {F_1} \subseteq Y_{F_2} = Y \cap {F_2}$,
where $Y = \{c_k | k \in J_N\}$. Similarly, we can prove that $Y_{F_2} \subseteq Y_{F_1}$. Therefore $Y_{F_1} = Y_{F_2}$.
Hence by Remark \ref{rem1}, $F_1 \cong F_2$. Thus $\phi$ is an injective map. 
\end{proof}
Note that, for every $\overrightarrow{G} \in \overrightarrow{\mathscr{D}}(n,l)$, there is a unique $G \in {\mathscr{D}}(n,l)$ which is obtained from $\overrightarrow{G}$ by removing direction of each edge $\overrightarrow{e} \in \overrightarrow{G}$. Also for every $G \in {\mathscr{D}}(n,l)$, there is a unique $\overrightarrow{G} \in \overrightarrow{\mathscr{D}}(n,l)$ with $V(\overrightarrow{G}) = V(G)$, and $(v_i, v_j) = \overrightarrow{e_k}$ is a directed edge in $\overrightarrow{G}$ whenever $k = (i-1)n-\binom{i}{2}+j-i$, $1 \leq i < j \leq n$. Thus there is a one-to-one correspondence between ${\mathscr{D}}(n,l)$ and $\overrightarrow{\mathscr{D}}(n,l)$. Thus by Theorem \ref{c6}, we get the following main result.
\begin{theorem} \label{MT}
For $n \geq 0$ and for fixed $l$, $\lfloor \frac{n+1}{2}\rfloor \leq l \leq \binom{n}{2}$, there is a one-to-one correspondence between $\mathscr{F}_n(l)$ and $\mathscr{D}(n,l)$, that is, $f(n, l) = d(n, l)$.
\end{theorem}
Thus from Theorem \ref{MT}, it clearly follows that the sequences $(f(n, l))$ and $(d(n, l))$ are equivalent.


\begin{thebibliography}{10}
\bibitem{bib1} A. N. Bhavale and B. N. Waphare, Basic retracts and counting of lattices, {\it Australas. J. Combin.}, \textbf{78}(2020), 73--99.
\bibitem{bib2} E. A. Bender, E. R. Canfield, and B. D. Mckay, The asymptotic number of labeled graphs with $n$ vertices, $q$ edges, and no isolated vertices, 
{\it J. Combin. Theory Ser. A}, \textbf{80}(1997), 124--150.
\bibitem{bib16} F. Brucker and A. G{\'e}ly, Crown-free lattices and their related graphs, {\it Order}, \textbf{3}(2011), 443--454.
\bibitem{bib9} E. N. Gilbert, Enumeration of labelled graphs, {\it Canad. J. Math.}, \textbf{8}(1956), 405--411.
\bibitem{bib15} G. Gr{\"a}tzer, {\it General Lattice Theory}, Birkh{\"a}user Verlag, 1978.
\bibitem{bib10} F. Harary, The number of linear, directed, rooted, and connected graphs, {\it Trans. Amer. Math. Soc.}, \textbf{78}(1955), 445--463.
\bibitem{bib11} F. Harary and E. Palmer, {\it Graphical Enumeration}, Academic Press, 1973.
\bibitem{bib4} I. Rival, Lattices with doubly irreducible elements, {\it Canad. Math. Bull.}, \textbf{17}(1974), 91--95.
\bibitem{bib20} A. Rosa, {\it On certain valuations of the vertices of a graph}, Theory of graphs (International Symposium, Rome, July 1996), Dunod Gordon and Breach, Paris (1967), 349--355.
\bibitem{bib3} N. J. A. Sloane, The Online Encyclopedia of Integer Sequences, https://oeis.org/, Concerned with sequence $A054548$.
\bibitem{bib13} R. Tauraso, Edge cover time for regular graphs, {\it J. Integer Seq.}, \textbf{11}(2008), Art.08.4.4.
\bibitem{bib5} N. K. Thakare, M. M. Pawar, and B. N. Waphare, A structure theorem for dismantlable lattices and enumeration, {\it Period. Math. Hungar.}, \textbf{45}(2002), 147--160.
\bibitem{bib14} D. B. West, {\it Introduction to Graph Theory}, Prentice Hall of India, 2002.
\end{thebibliography}
\end{document}